\documentclass{amsart}
\usepackage{amsmath}
\usepackage{amssymb}
\usepackage[dvips]{epsfig}
\usepackage{a4}		
\title{The Carmichael numbers up to $10^{17}$}
\author{Richard G.E. Pinch}
\address{2 Eldon Road, Cheltenham, Glos GL52 6TU, U.K.}
\email{rgep@chalcedon.demon.co.uk}
\date{5 April 2005}
\newcommand{\divides}{\vert}
\newcommand{\congruent}{\equiv}
\newcommand{\lcm}{\mathop{\rm lcm}}

\def\abs|#1|{{\left\vert{#1}\right\vert}}
\def\O(#1){\mbox{O}{\left({#1}\right)}}
\def\C(#1){C{\left({#1}\right)}}
\def\GF(#1){{\mathbb{F}_{#1}}}			
\def\gen<#1>{{\left\langle{#1}\right\rangle}}
\def\paren(#1){{\left({#1}\right)}}
\def\braces#1{{\left\lbrace{#1}\right\rbrace}}

\begin{document}

\maketitle

\begin{abstract}
We extend our previous computations to show that there are
585355 Carmichael numbers up to $10^{17}$.  As before,
the numbers were generated by a back-tracking search for
possible prime factorisations together with a
``large prime variation''.  We present further statistics on
the distribution of Carmichael numbers.
\end{abstract}

\suppressfloats

\section{Introduction}

A {\em Carmichael number} $N$ is a composite number $N$ with the property
that for every $b$ prime to $N$ we have $b^{N-1} \congruent 1 \mod N$.
It follows that a Carmichael number $N$
must be square-free, with at least three prime factors, and
that $p-1 \divides N-1$ for every prime $p$ dividing $N$:
conversely, any such $N$ must be a Carmichael number.

For background on Carmichael numbers and details of previous computations
we refer to our previous paper \cite{Pin:car15}: in that paper we described
the computation of the Carmichael numbers up to $10^{15}$ and presented
some statistics.  These computations have since been extended to 
$10^{16}$ \cite{Pin:car16} and now to 
$10^{17}$,
using similar techniques, and we present further statistics.

The complete list of Carmichael numbers up to $10^{17}$
is available from the author.

\section{Organisation of the search}

We used improved versions of strategies described in \cite{Pin:car15}.

The principal search was a depth-first back-tracking search over 
possible sequences of primes factors $p_1,\ldots,p_d$.  Put
$P_r = \prod_{i=1}^r p_i$, 
$Q_r = \prod_{i=r+1}^d p_i$ and
$L_r = \lcm\braces{p_i-1 : i=1,\ldots,r}$.
We find that $Q_r$ must satisfy the congruence 
$N = P_rQ_r \congruent 1 \bmod L_r$ and so
in particular $Q_d = p_d$ must satisfy a congruence modulo $L_{d-1}$:
further $p_d-1$ must be a factor of $P_{d-1}-1$.  We modified this to 
terminate the search early at some level $r$ if the modulus $L_r$ is
large enough to limit the possible values of $Q_r$, which 
may then be factorised directly.

We also employed the variant based on proposition 2 of \cite{Pin:car15}
which determines the finitely many possible pairs $(p_{d-1},p_d)$ from
$P_{d-2}$.  In practice this was useful only when $d=3$ allowing us to
determine the complete list of Carmichael numbers with three prime factors
up to $10^{18}$.  The count of these is 35586 (not 35585 as incorrectly
reported in our preliminary work).

Finally we employed a different search over large values of $p_d$, 
in the range $10^7 < p_d < 10^9$, using the property that 
$P_{d-1} \congruent 1 \bmod (p_d-1)$ and factorising the numbers in
this arithmetic progression less than $10^{18}/p_d$.  

\section{Statistics}

\begin{table}[phtb]
\begin{tabular}{||r|r||}
\hline
 $n$ & $C\left(10^n\right)$  \\
\hline
 3 &      1  \\
 4 &      7  \\
 5 &     16  \\
 6 &     43  \\
 7 &    105  \\
 8 &    255  \\
 9 &    646  \\
 10 &   1547  \\
 11 &   3605  \\
 12 &   8241  \\
 13 &  19279  \\
 14 &  44706  \\
 15 & 105212 \\
 16 & 246683 \\
 17 & 585355 \\
\hline
\end{tabular}
\caption{Distribution of Carmichael numbers up to $10^{17}$.}
\label{table0}
\end{table}

\begin{table}[phtb]
\begin{tabular}{||r|r|r|r|r|r|r|r|r|r|r||}
\hline
 $X$ &  3 &  4 &  5 &  6 &  7 &  8 &  9 &  10 & 11 &  total  \\
\hline
   3 &     1 &     0 &     0 &     0 &     0 &     0 &    0 &  0 & 0 &      1  \\
   4 &     7 &     0 &     0 &     0 &     0 &     0 &    0 &  0 & 0 &      7  \\
   5 &    12 &     4 &     0 &     0 &     0 &     0 &    0 &  0 & 0 &     16  \\
   6 &    23 &    19 &     1 &     0 &     0 &     0 &    0 &  0 & 0 &     43  \\
   7 &    47 &    55 &     3 &     0 &     0 &     0 &    0 &  0 & 0 &    105  \\
   8 &    84 &   144 &    27 &     0 &     0 &     0 &    0 &  0 & 0 &    255  \\
   9 &   172 &   314 &   146 &    14 &     0 &     0 &    0 &  0 & 0 &    646  \\
  10 &   335 &   619 &   492 &    99 &     2 &     0 &    0 &  0 & 0 &   1547  \\
  11 &   590 &  1179 &  1336 &   459 &    41 &     0 &    0 &  0 & 0 &   3605  \\
  12 &  1000 &  2102 &  3156 &  1714 &   262 &     7 &    0 &  0 & 0 &   8241  \\
  13 &  1858 &  3639 &  7082 &  5270 &  1340 &    89 &    1 &  0 & 0 &  19279  \\
  14 &  3284 &  6042 & 14938 & 14401 &  5359 &   655 &   27 &  0 & 0 &  44706  \\
  15 &  6083 &  9938 & 29282 & 36907 & 19210 &  3622 &  170 &  0 & 0 & 105212  \\
  16 & 10816 & 16202 & 55012 & 86696 & 60150 & 16348 & 1436 & 23 & 0 & 246683  \\
  17 & 19539 & 25758 &100707 &194306 &172234 & 63635 & 8835 &340 & 1 & 585355  \\
\hline
\end{tabular}
\caption{Values of $\C(X)$ and $\C(d,X)$ for $d \le 10$ and $X$ in
powers of 10 up to $10^{17}$.}
\label{table1}
\end{table}

We have shown that there are 585355 Carmichael numbers up to $10^{17}$,
all with at most 11 prime factors.  We let $\C(X)$ denote
the number of Carmichael numbers less than $X$ and $\C(d,X)$ denote the
number with exactly $d$ prime factors.  Table \ref{table0} gives the
values of $\C(X)$ and Table \ref{table1} the values
of $\C(d,X)$ for $X$ in powers of 10 up to $10^{17}$.

We have used the same methods to calculate the smallest Carmichael number
$S_d$ with $d$ prime factors for $d$ up to 34.  The results are given in
Table \ref{table2}.  Heuristics suggest that a good approximation for 
$\log S_d$ might be $2\log2 (d-1)\log(d-1)$.
There is some support for this from Table \ref{table2a}.

\begin{table}[phtb]
\begin{tabular}{||r|r||}
\hline
  $d$ &  $N$  \\ \hline
 3 & 561  \\ \hline
 4 & 41041  \\ \hline
 5 & 825265  \\ \hline
 6 & 321197185  \\ \hline
 7 & 5394826801  \\ \hline
 8 & 232250619601  \\ \hline
 9 & 9746347772161  \\ \hline
10 & 1436697831295441  \\ \hline
11 & 60977817398996785  \\ \hline
12 & 7156857700403137441  \\ \hline
13 & 1791562810662585767521  \\ \hline
14 & 87674969936234821377601  \\ \hline
15 & 6553130926752006031481761  \\ \hline
16 & 1590231231043178376951698401  \\ \hline
17 & 35237869211718889547310642241  \\ \hline
18 & 32809426840359564991177172754241  \\ \hline
19 & 2810864562635368426005268142616001  \\ \hline
20 & 349407515342287435050603204719587201  \\ \hline
21 & 125861887849639969847638681038680787361 \\ \hline
22 & 12758106140074522771498516740500829830401 \\ \hline
23 & 2333379336546216408131111533710540349903201 \\ \hline
24 & 294571791067375389885907239089503408618560001 \\ \hline
25 & 130912961974316767723865201454187955056178415601 \\ \hline
26 & 13513093081489380840188651246675032067011140079201 \\ \hline
27 & 7482895937713262392883306949172917048928068129206401 \\ \hline
28 & 1320340354477450170682291329830138947225695029536281601 \\ \hline
29 & 379382381447399527322618466130154668512652910714224209601 \\ \hline
30 & 70416887142533176417390411931483993124120785701395296424001 \\ \hline
31 & 2884167509593581480205474627684686008624483147814647841436801 \\ \hline
32 & 4754868377601046732119933839981363081972014948522510826417784001 \\ \hline
33 & 1334733877147062382486934807105197899496002201113849920496510541601 \\ \hline
34 & 260849323075371835669784094383812120359260783810157225730623388382401 \\ \hline
\end{tabular}
\caption{The smallest Carmichael numbers with $d$ prime factors for $d$ up to 34.}
\label{table2}
\end{table}

\begin{table}[phtb]
\begin{tabular}{||r|r||}
\hline
$d$ & $\frac{\log(S_d)}{2\log2 (d-1)\log(d-1)}$	\\	\hline
3 & 3.293621187 	\\
4 & 2.324869052 	\\
5 & 1.772215418 	\\
6 & 1.755823184 	\\
7 & 1.503593258 	\\
8 & 1.385943139 	\\
9 & 1.296862582 	\\
10 & 1.273113126 	\\
11 & 1.210794211 	\\
12 & 1.187291825 	\\
13 & 1.183842403 	\\
14 & 1.142841324 	\\
15 & 1.115637251 	\\
16 & 1.112255293 	\\
17 & 1.068846695 	\\
18 & 1.086833722 	\\
19 & 1.067861876 	\\
20 & 1.055266652 	\\
21 & 1.056211783 	\\
22 & 1.041906608 	\\
23 & 1.034833298 	\\
24 & 1.024202006 	\\
25 & 1.026042173 	\\
26 & 1.014073506 	\\
27 & 1.017122489 	\\
28 & 1.010169037 	\\
29 & 1.007224804 	\\
30 & 1.000945160 	\\
31 & 0.984182030 	\\
32 & 0.993535534 	\\
33 & 0.990337138 	\\
34 & 0.984854273 	\\ \hline
\end{tabular}
\caption{The smallest Carmichael number with $d$ prime factors $S_d$ 
compared to $2\log2 (d-1)\log(d-1)$.}
\label{table2a}
\end{table}

\begin{figure}
\epsfig{file=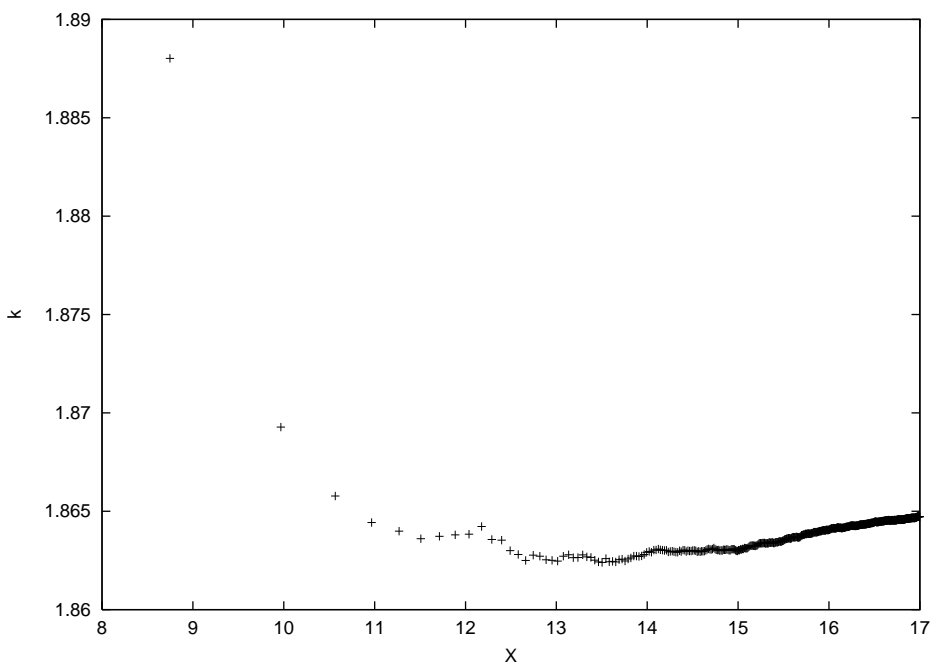,width=\textwidth}
\vspace{-1truecm}
\caption{$k(X)$ versus $X$ (expressed in powers of 10).}
\label{figure1}
\end{figure}

\begin{table}[phtb]
\begin{tabular}{||r|r|c||}
\hline
   $n$ &  $k\left(10^n\right)$ &  $\C(10^n)/\C(10^{n-1})$ \\
\hline
    3 &  2.93319 &        \\
    4 &  2.19547 &  7.000 \\
    5 &  2.07632 &  2.286 \\
    6 &  1.97946 &  2.688 \\
    7 &  1.93388 &  2.441 \\
    8 &  1.90495 &  2.429 \\
    9 &  1.87989 &  2.533 \\
   10 &  1.86870 &  2.396 \\
   11 &  1.86421 &  2.330 \\
   12 &  1.86377 &  2.286 \\
   13 &  1.86240 &  2.339 \\
   14 &  1.86293 &  2.319 \\
   15 &  1.86301 &  2.353 \\
   16 &  1.86406 &  2.335 \\
   17 &  1.86472 &  2.373 \\
\hline
\end{tabular}
\caption{The function $k(X)$ and growth of $\C(X)$ for $X = 10^n$,
$n \le 17$.}
\label{table3}
\end{table}

In Table \ref{table3} and Figure \ref{figure1} we tabulate the function
$k(X)$, defined by Pomerance, Selfridge and Wagstaff \cite{PSW:pseudo} by
$$
\C(X) = X \exp\left(-k(X) \frac{\log X \log\log\log X}{\log\log X}\right) .
$$
They proved that $\liminf k \ge 1$
and suggested that $\limsup k$ might be $2$, although they also observed
that within the range of their tables $k(X)$ is decreasing:
Pomerance \cite{Pom:distpsp},\cite{Pom:2meth}
gave a heuristic argument suggesting that $\lim k = 1$.
The decrease in $k$ is reversed between $10^{13}$ and $10^{14}$: see
Figure \ref{figure1}.  We find no clear support from our computations
for any conjecture on a limiting value of $k$.

In Table \ref{table3} we also give the ratios $\C(10^n) / \C(10^{n-1})$
investigated by Swift \cite{Swi:car}.
Swift's ratio, again initially decreasing, also increases again before
$10^{15}$.

\begin{figure}
\epsfig{file=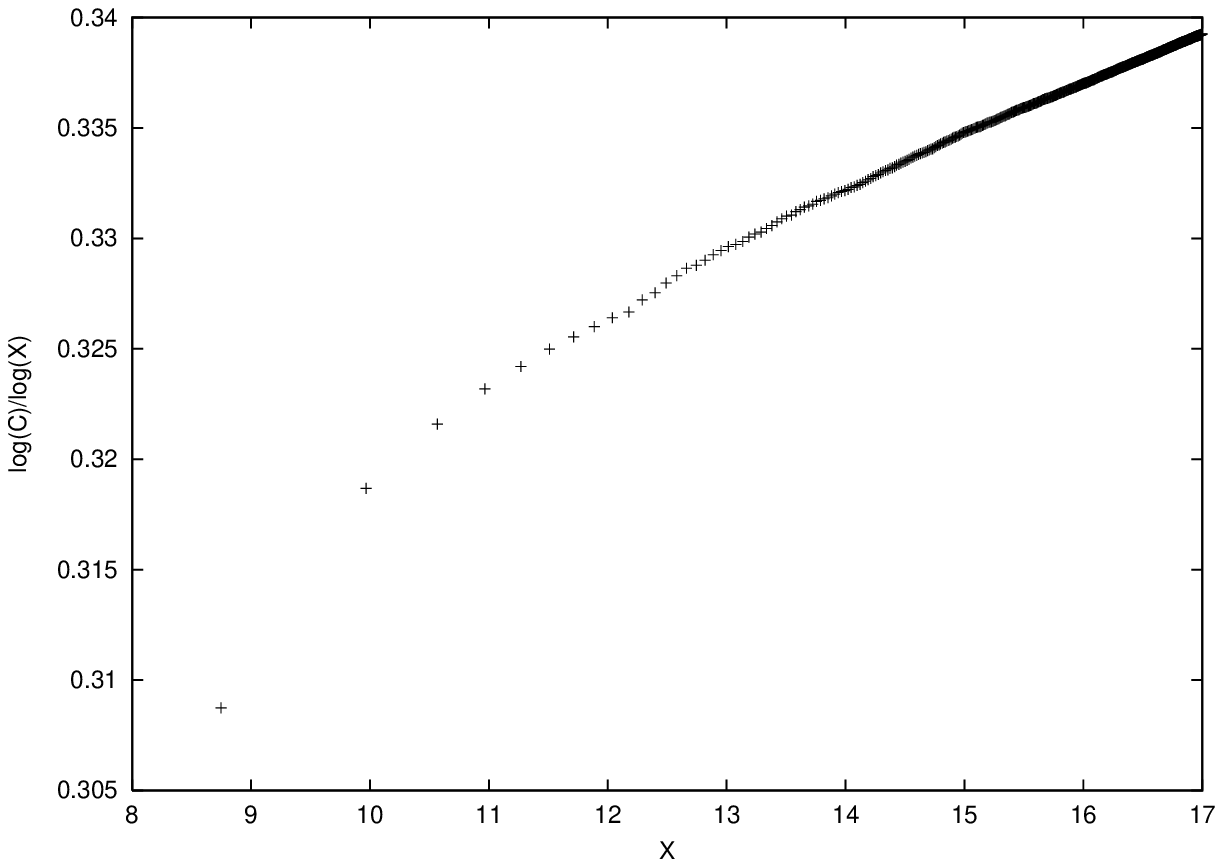,width=\textwidth}
\vspace{-1truecm}
\caption{$C(X)$ as a power of $X$ (expressed in powers of 10).}
\label{figure2}
\end{figure}

\begin{table}[phtb]
\begin{tabular}{||l|c||}
\hline
 $n$ &  $\log C\left(10^n\right) / (n \log 10)$ \\
\hline
 4  &  .21127  \\
 5  &  .24082  \\
 6  &  .27224  \\
 7  &  .28874  \\
 8  &  .30082  \\
 9  &  .31225  \\
 {10} &  .31895  \\
 {11} &  .32336  \\
 {12} &  .32633  \\
 {13} &  .32962  \\
 {14} &  .33217  \\
 {15} &  .33480  \\
 {16} &  .33700  \\
 {17} &  .33926  \\
\hline
\end{tabular}
\caption{$\C(X)$ as a power of $X$.}
\label{table00}
\end{table}

In Table \ref{table00} and Figure \ref{figure2} we see that within the
range of our computations $C(X)$ is a slowly growing power of $X$:
about $X^{0.339}$ for $X$ around $10^{17}$.

In Table \ref{table4} we give the number of Carmichael numbers in each
class modulo $m$ for $m = 5$, 7, 11 and 12.

In Tables \ref{table5} and \ref{table6} we give the number of Carmichael
numbers divisible by primes $p$ up to 97.  In Table \ref{table5} we count
all Carmichael numbers divisible by $p$: in Table \ref{table6} we count
only those for which $p$ is the smallest prime factor.

The largest prime factor of a Carmichael number up
to $10^{17}$ is $p=223401361$, 
dividing $ n = 99816335969903281 = 17 \cdot 31 \cdot 379 \cdot 2237 \cdot p$.
We have $p = n^{0.4911}$: indeed $n = p(2p-1)$ and $n-1 = (p-1)(2p+1)$.

The largest prime to occur as the smallest prime factor of a Carmichael
number in this range is 380251, dividing
$$
90256390764228001 = 380251 \cdot 410761 \cdot 577981 .
$$
We note that this number is of the form $(25k+1)(27k+1)(38k+1)$
with $k = 15210$.

We define the {\em index} of a Carmichael number $N$ to be the integer
$i(N) = (N-1)/\lambda(N)$ where $\lambda$ is the Carmichael function, the
exponent of the multiplicative group.  A result of Somer \cite{Som:Ferdpsp}
implies that $i(N) \rightarrow \infty$ as $N$ runs through Carmichael
numbers.  In Table \ref{table7} we list the Carmichael numbers known to
have index less than 100: the list is complete for numbers up to $10^{17}$.

We define the {\em Lehmer index} of a number $N$ to be $\ell(N) = (N-1)/\phi(N)$.
Lehmer's totient problem asks whether there is a number of integral
Lehmer index greater than 1 (that is, whether there is a composite $N$
such that $\phi(N)$ divides $N-1$).  Such a number would have to be a
Carmichael number.  See Guy \cite{Guy:upint} \S B37 for further details.

The Carmichael numbers up to $10^{17}$ with Lehmer index $\ge 2$ are listed 
in Table \ref{table8}.

\begin{table}[phtb]
\begin{tabular}{||r|r|r|r|r|r|r|r|r|r||}
\hline
 $m$ & $c$ & $25.10^9$ & $10^{11}$ & $10^{12}$ & $10^{13}$ & $10^{14}$ & $10^{15}$ & $10^{16}$ & $10^{17}$ \\
\hline
 5  &  0 & 203 &  312 &  627 & 1330 &  2773 &  5814 &  12200 &  25523 \\
    &  1 &1652 & 2785 & 6575 &15755 & 37467 & 90167 & 215713 & 520990 \\
    &  2 &  82 &  154 &  327 &  702 &  1484 &  3048 &   6094 &  12912 \\
    &  3 & 102 &  172 &  344 &  725 &  1463 &  3059 &   6285 &  12845 \\
    &  4 & 124 &  182 &  368 &  767 &  1519 &  3124 &   6391 &  13085 \\
\hline
 7  &  0 & 401 &  634 & 1334 & 2774 &  5891 & 12691 &  27550 &  60602 \\
    &  1 &1096 & 1885 & 4613 &11447 & 28001 & 69131 & 168856 & 414999 \\
    &  2 & 105 &  186 &  432 &  967 &  2109 &  4599 &  10011 &  21961 \\
    &  3 & 152 &  232 &  496 & 1055 &  2178 &  4707 &  10039 &  21918 \\
    &  4 & 129 &  211 &  450 &  985 &  2122 &  4592 &   9944 &  21832 \\
    &  5 & 138 &  222 &  454 & 1033 &  2224 &  4777 &  10125 &  21911 \\
    &  6 & 142 &  235 &  462 & 1018 &  2181 &  4715 &  10158 &  22132 \\
\hline
11  &  0 & 335 &  547 & 1324 & 3006 &  7032 & 16563 &  38576 &  90731 \\
    &  1 & 640 & 1131 & 2770 & 6786 & 16548 & 40891 & 100071 & 248840 \\
    &  2 & 139 &  217 &  473 & 1068 &  2361 &  5338 &  12054 &  27268 \\
    &  3 & 142 &  220 &  457 & 1045 &  2348 &  5319 &  12186 &  27483 \\
    &  4 & 104 &  187 &  442 & 1026 &  2317 &  5261 &  11917 &  27276 \\
    &  5 & 152 &  243 &  466 & 1066 &  2370 &  5316 &  12194 &  27670 \\
    &  6 & 116 &  198 &  440 & 1061 &  2400 &  5384 &  12155 &  27362 \\
    &  7 & 122 &  195 &  458 & 1023 &  2223 &  5165 &  11853 &  26960 \\
    &  8 & 129 &  222 &  475 & 1107 &  2450 &  5449 &  12012 &  27639 \\
    &  9 & 131 &  218 &  465 & 1042 &  2285 &  5179 &  11835 &  27158 \\
    & 10 & 153 &  227 &  471 & 1049 &  2372 &  5347 &  11830 &  26968 \\
\hline
12  &  1 &2071 &3462 &7969 &18761 &43760 &103428 &243382 & 579215 \\
    &  3 &   0 &   0 &   1 &    2 &    2 &     5 &     5 &      7 \\
    &  5 &  20 &  32 &  64 &  124 &  228 &   448 &   805 &   1470 \\
    &  7 &  47 &  75 & 147 &  289 &  547 &  1027 &  1906 &   3589 \\
    &  9 &  25 &  36 &  60 &  103 &  165 &   294 &   560 &   1018 \\
    & 11 &   0 &   0 &   0 &    0 &    4 &    10 &    25 &     56 \\
\hline
\end{tabular}
\caption{The number of Carmichael numbers in each class modulo $m$
for $m = 5$, 7, 11 and 12.}
\label{table4}
\end{table}

\begin{table}[phtb]
\begin{tabular}{||r|r|r|r|r|r|r|r|r|r||}
\hline
$p$ & $25.10^9$ & $10^{11}$ & $10^{12}$ & $10^{13}$ & $10^{14}$ & $10^{15}$ & $10^{16}$ & $10^{17}$ \\
\hline
 3 &         25 &       36 &       61 &      105 &      167 &      299 &      565 &      1025 \\
 5 &        203 &      312 &      627 &     1330 &     2773 &     5814 &    12200 &     25523 \\
 7 &        401 &      634 &     1334 &     2774 &     5891 &    12691 &    27550 &     60602 \\
11 &        335 &      547 &     1324 &     3006 &     7032 &    16563 &    38576 &     90731 \\
13 &        483 &      807 &     1784 &     3998 &     9045 &    20758 &    47785 &    110320 \\
17 &        293 &      489 &     1182 &     2817 &     6640 &    16019 &    38302 &     91783 \\
19 &        372 &      608 &     1355 &     3345 &     7797 &    18638 &    44389 &    106273 \\
23 &        113 &      207 &      507 &     1282 &     3135 &     7716 &    18867 &     46612 \\
29 &        194 &      336 &      832 &     2094 &     5158 &    12721 &    31110 &     76647 \\
31 &        335 &      571 &     1320 &     3086 &     7270 &    17382 &    41440 &     99845 \\
37 &        320 &      535 &     1270 &     2926 &     6826 &    16220 &    38647 &     92744 \\
41 &        227 &      390 &     1001 &     2418 &     5896 &    14344 &    34759 &     84977 \\
43 &        184 &      296 &      772 &     1920 &     4663 &    11594 &    28650 &     70904 \\
47 &         53 &       80 &      199 &      492 &     1223 &     2873 &     6810 &     17002 \\
53 &         92 &      160 &      351 &      813 &     2041 &     5143 &    12256 &     30560 \\
59 &         26 &       41 &       92 &      262 &      644 &     1611 &     3959 &      9613 \\
61 &        269 &      453 &     1075 &     2542 &     6047 &    14429 &    34503 &     83268 \\
67 &        110 &      178 &      407 &     1063 &     2540 &     6306 &    15295 &     37521 \\
71 &        104 &      194 &      521 &     1320 &     3351 &     8546 &    21485 &     53857 \\
73 &        198 &      348 &      849 &     2145 &     4925 &    11929 &    29072 &     71137 \\
79 &         64 &      107 &      247 &      686 &     1728 &     4318 &    10693 &     26766 \\
83 &         14 &       24 &       56 &      137 &      340 &      838 &     1929 &      4640 \\
89 &         68 &      131 &      320 &      788 &     1951 &     4981 &    12178 &     30737 \\
97 &        123 &      193 &      495 &     1277 &     3123 &     7594 &    18706 &     45542 \\
\hline
\end{tabular}
\caption{Primes occurring in Carmichael numbers.}
\label{table5}
\end{table}

\begin{table}[phtb]
\begin{tabular}{||r|r|r|r|r|r|r|r|r||}
\hline
$p$ & $25.10^9$ & $10^{11}$ & $10^{12}$ & $10^{13}$ & $10^{14}$ & $10^{15}$ & $10^{16}$ & $10^{17}$ \\
\hline
 3 &        25 &      36 &      61 &     105 &     167 &     299 &     565 &    1025 \\
 5 &       202 &     309 &     624 &    1325 &    2765 &    5797 &   12175 &   25481 \\
 7 &       364 &     579 &    1218 &    2557 &    5461 &   11874 &   25915 &   57459 \\
11 &       263 &     428 &    1071 &    2509 &    5979 &   14397 &   33893 &   80745 \\
13 &       237 &     431 &    1058 &    2462 &    5699 &   13514 &   32025 &   76256 \\
17 &       117 &     206 &     496 &    1318 &    3244 &    8114 &   20206 &   50170 \\
19 &       152 &     244 &     532 &    1401 &    3358 &    8141 &   20020 &   49413 \\
23 &        37 &      78 &     207 &     535 &    1360 &    3317 &    8195 &   20803 \\
29 &        55 &     103 &     284 &     729 &    1822 &    4659 &   11577 &   29149 \\
31 &       101 &     168 &     390 &     876 &    2116 &    5153 &   12575 &   30667 \\
37 &        60 &      95 &     219 &     551 &    1401 &    3418 &    8594 &   21382 \\
41 &        35 &      68 &     171 &     414 &    1092 &    2736 &    6788 &   17275 \\
43 &        35 &      65 &     168 &     403 &     943 &    2308 &    5520 &   13636 \\
47 &        14 &      16 &      36 &      81 &     195 &     459 &    1135 &    2854 \\
53 &        19 &      30 &      55 &     147 &     363 &     973 &    2327 &    5842 \\
59 &         2 &       4 &      11 &      43 &     100 &     272 &     618 &    1542 \\
61 &        34 &      58 &     148 &     364 &     851 &    1978 &    4722 &   11278 \\
67 &         8 &      18 &      50 &     123 &     317 &     815 &    1950 &    4843 \\
71 &        15 &      25 &      66 &     161 &     389 &     979 &    2480 &    6178 \\
73 &        14 &      28 &      68 &     175 &     406 &    1015 &    2508 &    6277 \\
79 &         4 &      10 &      17 &      66 &     175 &     467 &    1163 &    2873 \\
83 &         1 &       1 &       4 &       8 &      39 &      79 &     175 &     457 \\
89 &        10 &      16 &      23 &      55 &     148 &     409 &    1003 &    2523 \\
97 &        10 &      20 &      50 &     106 &     261 &     606 &    1413 &    3445 \\
\hline
\end{tabular}
\caption{Primes occurring as least prime factor in Carmichael numbers.}
\label{table6}
\end{table}

\clearpage

\begin{table}[phtb]
\begin{tabular}{||r|r|l||}
\hline
$i$ & $N$ & factors \\
\hline
5 & 6601 & $ 7 \cdot 23 \cdot 41 $	\\
7 & 561 & $ 3 \cdot 11 \cdot 17 $	\\
18 & 42018333841 & $ 11 \cdot 47 \cdot 1049 \cdot 77477 $	\\
18 & 55462177 & $ 17 \cdot 23 \cdot 83 \cdot 1709 $	\\
18 & 8885251441 & $ 11 \cdot 47 \cdot 1109 \cdot 15497 $	\\
21 & 10585 & $ 5 \cdot 29 \cdot 73 $	\\
22 & 2465 & $ 5 \cdot 17 \cdot 29 $	\\
23 & 1105 & $ 5 \cdot 13 \cdot 17 $	\\
25 & 11921001 & $ 3 \cdot 29 \cdot 263 \cdot 521 $	\\
31 & 62745 & $ 3 \cdot 5 \cdot 47 \cdot 89 $	\\
37 & 11972017 & $ 43 \cdot 433 \cdot 643 $	\\
37 & 67902031 & $ 43 \cdot 271 \cdot 5827 $	\\
39 & 334153 & $ 19 \cdot 43 \cdot 409 $	\\
43 & 52633 & $ 7 \cdot 73 \cdot 103 $	\\
44 & 15841 & $ 7 \cdot 31 \cdot 73 $	\\
45 & 8911 & $ 7 \cdot 19 \cdot 67 $	\\
47 & 2821 & $ 7 \cdot 13 \cdot 31 $	\\
48 & 1729 & $ 7 \cdot 13 \cdot 19 $	\\
49 & 1208361237478669 & $ 53 \cdot 653 \cdot 26479 \cdot 1318579 $	\\
50 & 4199932801 & $ 29 \cdot 499 \cdot 503 \cdot 577 $	\\
52 & 206955841 & $ 17 \cdot 71 \cdot 277 \cdot 619 $	\\
53 & 1271325841 & $ 17 \cdot 31 \cdot 179 \cdot 13477 $	\\
54 & 4169867689 & $ 13 \cdot 29 \cdot 383 \cdot 28879 $	\\
55 & 271794601 & $ 13 \cdot 19 \cdot 743 \cdot 1481 $	\\
60 & 6840001 & $ 7 \cdot 17 \cdot 229 \cdot 251 $	\\
61 & 1962804565 & $ 5 \cdot 103 \cdot 149 \cdot 25579 $	\\
67 & 11985924995083901 & $ 29 \cdot 101 \cdot 1427 \cdot 16349 \cdot 175403 $	\\
67 & 410041 & $ 41 \cdot 73 \cdot 137 $	\\
70 & 162401 & $ 17 \cdot 41 \cdot 233 $	\\
76 & 4752717761 & $ 11 \cdot 17 \cdot 107 \cdot 173 \cdot 1373 $	\\
81 & 35575075809505 & $ 5 \cdot 197 \cdot 223 \cdot 353 \cdot 458807 $	\\
82 & 1496405933740345 & $ 5 \cdot 47 \cdot 317 \cdot 40253 \cdot 499027 $	\\
83 & 142159958924185 & $ 5 \cdot 37 \cdot 107 \cdot 58379 \cdot 123017 $	\\
83 & 158664761899885 & $ 5 \cdot 37 \cdot 107 \cdot 53987 \cdot 148469 $	\\
83 & 204370370140285 & $ 5 \cdot 37 \cdot 107 \cdot 48677 \cdot 212099 $	\\
83 & 24831908105124205 & $ 5 \cdot 29 \cdot 719 \cdot 3023 \cdot 78790717 $	\\
85 & 171189355538562901 & $ 23 \cdot 71 \cdot 983 \cdot 1031 \cdot 103437589 $	\\
90 & 3778118040573702001 & $ 11 \cdot 47 \cdot 1051 \cdot 67967 \cdot 102302009 $ 	\\
92 & 520178982961 & $ 11 \cdot 29 \cdot 131 \cdot 607 \cdot 20507 $	\\
94 & 12782849065 & $ 5 \cdot 7 \cdot 269 \cdot 317 \cdot 4283 $	\\
95 & 8956911601 & $ 11 \cdot 17 \cdot 127 \cdot 131 \cdot 2879 $	\\
97 & 5472940991761 & $ 199 \cdot 241 \cdot 863 \cdot 132233 $	\\
97 & 721574219707441 & $ 167 \cdot 241 \cdot 5039 \cdot 3557977 $	\\
97 & 83565865434172201 & $ 103 \cdot 1993 \cdot 9551 \cdot 42622169 $	\\
97 & 9729822470631481 & $ 127 \cdot 409 \cdot 110681 \cdot 1692407 $	\\
99 & 438253965870337 & $ 43 \cdot 139 \cdot 49409 \cdot 1484009 $	\\
\hline
\end{tabular}
\caption{Carmichael numbers with small index}
\label{table7}
\end{table}
	
\begin{table}[phtb]
\begin{tabular}{||r|r|l||}
\hline
$\ell$ & $N$ & factors \\
\hline
2.00163  &  3852971941960065
&  $3 \cdot 5 \cdot 23 \cdot 89 \cdot 113 \cdot 1409 \cdot 788129$	\\
2.00388  &  655510549443465
&  $3 \cdot 5 \cdot 23 \cdot 53 \cdot 389 \cdot 2663 \cdot 34607$	\\
2.00837    &    13462627333098945
&  $3 \cdot 5 \cdot 23 \cdot 53 \cdot 197 \cdot 8009 \cdot 466649$	\\
2.01001    &    26708253318968145
&  $3 \cdot 5 \cdot 17 \cdot 113 \cdot 57839 \cdot 16025297$	\\
2.03067    &    269040992399565
&  $3 \cdot 5 \cdot 23 \cdot 29 \cdot 4637 \cdot 5799149$	\\
2.03207    &    158353658932305
&  $3 \cdot 5 \cdot 17 \cdot 89 \cdot 149 \cdot 563 \cdot 83177$	\\
2.03614    &    1817671359979245
&  $3 \cdot 5 \cdot 23 \cdot 29 \cdot 359 \cdot 11027 \cdot 45893$	\\
2.07138    &    16057190782234785
&  $3 \cdot 5 \cdot 17 \cdot 29 \cdot 269 \cdot 6089 \cdot 1325663$	\\
2.07294    &    75131642415974145
&  $3 \cdot 5 \cdot 23 \cdot 29 \cdot 53 \cdot 617 \cdot 9857 \cdot 23297$	\\
2.08129    &    881715504450705
&  $3 \cdot 5 \cdot 17 \cdot 47 \cdot 89 \cdot 113 \cdot 503 \cdot 14543$	\\
2.08381    &    31454143858820145
&  $3 \cdot 5 \cdot 17 \cdot 23 \cdot 2129 \cdot 39293 \cdot 64109$	\\
2.08894    &    6128613921672705
&  $3 \cdot 5 \cdot 17 \cdot 23 \cdot 353 \cdot 7673 \cdot 385793$	\\
2.09841    &    12301576752408945
&  $3 \cdot 5 \cdot 23 \cdot 29 \cdot 53 \cdot 113 \cdot 197 \cdot 1042133$	\\
2.11432    &    1886616373665
&  $3 \cdot 5 \cdot 17 \cdot 23 \cdot 83 \cdot 353 \cdot 10979$	\\
2.11493    &    3193231538989185
&  $3 \cdot 5 \cdot 17 \cdot 23 \cdot 113 \cdot 167 \cdot 2927 \cdot 9857$	\\
2.13936    &    11947816523586945
&  $3 \cdot 5 \cdot 17 \cdot 23 \cdot 89 \cdot 113 \cdot 233 \cdot 617 \cdot 1409$	\\
\hline
\end{tabular}
\caption{Carmichael numbers with Lehmer index greater than 2}
\label{table8}
\end{table}
	

\nocite{Mol:ntapps}\nocite{KL:Laval87}

\providecommand{\bysame}{\leavevmode\hbox to3em{\hrulefill}\thinspace}
\providecommand{\MR}{\relax\ifhmode\unskip\space\fi MR }
\providecommand{\MRhref}[2]{%
  \href{http://www.ams.org/mathscinet-getitem?mr=#1}{#2}
}
\providecommand{\href}[2]{#2}

\end{document}